\documentclass[a4paper,12pt]{article}

\usepackage[T1]{fontenc}
\usepackage[utf8]{inputenc}
\usepackage[english]{babel}
\usepackage[inline,shortlabels]{enumitem}
\usepackage{amsmath}
\usepackage{amssymb}
\usepackage{amsthm}
\usepackage[noadjust]{cite}
\usepackage[a4paper,total={6in,9in}]{geometry}
\usepackage[colorlinks,citecolor=blue]{hyperref}
\usepackage[capitalize,noabbrev,nameinlink]{cleveref}
\usepackage{xcolor}
\usepackage{longtable}
\usepackage{tikz,pgfplots,pgfplotstable}
\usetikzlibrary{patterns}
\usepackage{algorithmic}
\usepackage{marginnote}
\usepackage[labelformat=simple]{subcaption}

\setlist{font=\normalfont,topsep=1ex,parsep=0ex}

\setlist[enumerate]{label=(\alph*)}

\numberwithin{equation}{section}
\numberwithin{table}{section}    
\numberwithin{figure}{section}

\crefname{figure}{Figure}{Figures}
\crefname{table}{Table}{Tables}
\crefname{assumption}{Assumption}{Assumptions}
\Crefname{ALC@unique}{Step}{Steps}

\newcommand{\R}{\mathbb{R}}

\newcommand\norm[1]{\left\Vert#1\right\Vert}
\newcommand{\N}{\mathbb{N}}

\DeclareMathOperator{\dist}{dist}

\newcommand{\dom}{\operatorname{dom}}

\newtheoremstyle{bolddef}{}{}{\normalfont}{}{\bfseries}{.}{ }{\thmname{#1}\thmnumber{ #2}\thmnote{ (#3)}}

\newtheoremstyle{boldplain}{}{}{\itshape}{}{\bfseries}{.}{ }{\thmname{#1}\thmnumber{ #2}\thmnote{ (#3)}}

\theoremstyle{bolddef}
\newtheorem{definition}{Definition}[section]
\newtheorem{algorithm}[definition]{Algorithm}
\newtheorem{assumption}[definition]{Assumption}
\newtheorem{example}[definition]{Example}

\theoremstyle{boldplain}
\newtheorem{lemma}[definition]{Lemma}
\newtheorem{theorem}[definition]{Theorem}
\newtheorem{proposition}[definition]{Proposition}

\newlength\figureheight
\newlength\figurewidth

\pgfplotsset{width=7cm,compat=1.3}

\definecolor{todocolor}{rgb}{1.0,0.0,0.0}

\newcommand\email[1]{\href{mailto:#1}{\texttt{#1}}}

\newcommand{\orcid}[1]{ORCID: \href{https://orcid.org/#1}{#1}}

\newcommand{\mscLink}[1]{\href{http://www.ams.org/mathscinet/msc/msc2020.html?t=#1}{#1}}

\begin{document}

\title{
	\bfseries\scshape 
	Convergence Analysis of the Proximal Gradient Method 
	in the Presence of the Kurdyka--{\L}ojasiewicz Property
	without Global Lipschitz Assumptions
	}

\date{\today}

\author{Xiaoxi Jia%
	\thanks{%
		University of Würzburg,
		Institute of Mathematics,
		97074 Würzburg,
		Germany,
		\email{xiaoxi.jia@mathematik.uni-wuerzburg.de}
	}
	\hspace*{-4mm} \and \hspace*{-4mm}
	Christian Kanzow%
	\thanks{%
		University of Würzburg,
		Institute of Mathematics,
		97074 Würzburg,
		Germany,
		\email{kanzow@mathematik.uni-wuerzburg.de},
		\orcid{0000-0003-2897-2509}
	}
	\hspace*{-4mm} \and \hspace*{-4mm}
	Patrick Mehlitz%
	\thanks{%
	Brandenburgische Technische Universität Cottbus-Senftenberg,
		Institute of Mathematics,
		03046 Cottbus,
		Germany,
		\email{mehlitz@b-tu.de},
		\orcid{0000-0002-9355-850X}
	} 
}

\maketitle

{
\small\textbf{\abstractname.}
We consider a composite optimization problem where the 
sum of a continuously differentiable and a merely lower semicontinuous function 
has to be minimized. 
The proximal gradient algorithm is the classical method for solving such a
problem numerically. The corresponding global convergence
and local rate-of-convergence theory typically assumes, besides some
technical conditions, that the smooth function has
a globally Lipschitz continuous gradient and that the objective
function satisfies the Kurdyka--{\L}ojasiewicz property. Though
this global Lipschitz assumption is satisfied in several applications
where the objective function is, e.g., quadratic, this requirement
is very restrictive in the non-quadratic case. Some recent
contributions therefore try to overcome this global Lipschitz 
condition by replacing it with a local one, but,
to the best of our knowledge, they still require some extra
condition in order to obtain the desired global and rate-of-convergence
results. The aim of this paper is to show that the local Lipschitz
assumption together with the Kurdyka--{\L}ojasiewicz property is
sufficient to recover these convergence results.
\par\addvspace{\baselineskip}
}

{
\par\addvspace{\baselineskip}
}

{
\par\addvspace{\baselineskip}
}

\section{Introduction}\label{Sec:Intro}

In this paper, we are concerned with problems from \emph{composite optimization}
where the sum of a continuously differentiable function $f$ and a merely
lower semicontinuous function $\phi$ has to be minimized.
Problems of this type appear quite frequently in 
many practically relevant areas like, e.g.,
machine learning, 
data compression, matrix completion, and image processing, see
\cite{BianChen2015,BrucksteinDonohoElad2009,Chartrand2007,DiLorenzoLiuzziRinaldiSchoenSciandrone2012,LiuDaiMa2015,MarjanovicSolo2012},
where, typically, $f$ models a tracking-type term while $\phi$ is used to promote sparse structures in the solutions.

For an algorithmic treatment of such problems, one can
exploit the composite form,
i.e., differentiability of $f$ on the one hand and additional structural properties 
of the function $\phi$ on the other hand (typically, the nonsmoothness encapsulated 
within $\phi$ is of specific type in all aforementioned applications).
More precisely, the so-called proximal mapping of the function $\phi$ has to be available,
which is typically the case in the aforementioned practically relevant scenarios. 
The idea behind the definition of proximal mappings is to interrelate the search for
minimizers (or at least stationary points) with a fixed-point problem, and to apply
a fixed-point iteration to the proximal mapping in order to tackle the minimization of
the underlying function.
Combining the available oracles for $f$ and $\phi$ in order to construct an algorithm to
minimize $f+\phi$ led to the development of so-called proximal gradient methods which date back
to \cite{FukushimaMine1981}. 
It is worth noting that proximal gradient algorithms can be interpreted as so-called
forward-backward splitting methods which are far older, see
\cite{Bruck1977,BruckReich1977,Passty1979,Rockafellar1976} for their origins and \cite{BauschkeCombettes2017} for a modern view. 
Popular instances of proximal gradient methods are the 
iterative shrinkage/threshold algorithm (ISTA) and its
accelerated version (FISTA = fast ISTA), see \cite{BeckTeboulle2009}, where $\phi$
has to be convex. The monograph \cite{Beck2017} presents a nice overview of existing
results addressing proximal gradient methods where the nonsmooth part enjoys convexity.

It has been pointed out in the seminal works 
\cite{AttouchBolteSvaiter2013,BolteSabachTeboulle2014}
that the convergence theory for proximal gradient methods can be extended to situations where
the nonsmooth part $\phi$ is merely lower semicontinuous and not necessarily convex.
In both aforementioned papers, the analysis, which covers both (global) convergence 
and rate-of-convergence results, requires a so-called \emph{descent lemma} 
as well as the celebrated \emph{Kurdyka--{\L}ojasiewicz property}, 
originating from \cite{Kurdyka1998,Lojasiewicz1963,Lojasiewicz1965}.
The majority of available convergence results regarding proximal gradient methods
seems to indicate that the price we have to pay for allowing $\phi$ to be nonsmooth
is that the gradient $\nabla f$ of the smooth part has to be globally Lipschitz continuous.
This requirement, which holds naturally when $f$ is a (convex) quadratic function
(as indicated above, this happens to be the case in many standard applications 
from image processing and data science), turns out to be rather restrictive in the non-quadratic
situation which also is of practical interest, see \cref{ex:alm,ex:dual_prox_grad} below.

Let us review some contributions 
where the authors try to get rid of this global Lipschitz assumption.
First, we would like to mention \cite{BauschkeBolteTeboulle2017}
where composite optimization problems with convex functions $f$ and $\phi$ are
considered without postulating global Lipschitzness of $\nabla f$.
It is shown that local Lipschitz continuity of $\nabla f$ is enough to obtain
rate-of-convergence results for the iterates generated
by a Bregman-type proximal gradient method. 
However, the authors of \cite{BauschkeBolteTeboulle2017} require the additional 
assumption that there is a constant $ L > 0 $ such that $ L h - f $ is
convex, where $ h $ is a convex function which defines the
Bregman distance (let us mention that $ h $ equals the squared 
Euclidean norm in our setting). 
This convexity-type condition is satisfied in a couple of practically relevant situations. 
The approach of \cite{BauschkeBolteTeboulle2017} was generalized to the nonconvex setting in 
\cite{BolteSabachTeboulleVaisbourd2018} using, once again,
a local Lipschitz assumption on $\nabla f$, as well as the slightly
stronger assumption (in order to deal with the nonconvexity)
that there exist a constant $ L > 0 $ 
and a convex function $ h $ such
that both $ L h - f $ and $ L h + f $ are convex. 
Let us emphasize that this constant $ L $ plays a central role in the design of the
corresponding proximal-type methods. More precisely, it is used explicitly
for the determination of the stepsizes.
In the recent paper \cite{CohenHallakTeboulle2022}, global convergence
results are proven under a local Lipschitz assumption on $\nabla f$ 
(without postulating any of the convexity-type conditions from above), but the authors assume
(a priori) boundedness of iterates and stepsizes.
Let us also mention some related works which do not address proximal gradient algorithms.
In \cite{KhanhMordukhovichTran2022,NollRondepierre2013}, the authors are concerned with
(inexact) descent methods for differentiable functions without a Lipschitzian gradient
and also investigate situations where the aforementioned Kurdyka--{\L}ojasiewicz property
is present. The paper \cite{Pauwels2016} studies the (convex-)constrained minimization of the composition
of a convex and a twice continuously differentiable function whose gradient is not assumed to be
globally Lipschitzian, based on a (nonsmooth) Gau\ss--Newton method. 
Using the Kurdyka--{\L}ojasiewicz property, convergence of the whole sequence of iterates is shown.

The present paper is based on \cite{KanzowMehlitz2022} where the
authors show global convergence results for proximal gradient
methods in the sense that every accumulation point is shown to
be a suitable stationary point of the composite optimization
problem. 
The analysis in \cite{KanzowMehlitz2022} is based on the local Lipschitz
continuity of $\nabla f$, and does not require the iterates to be bounded.
	Related results under similar assumptions for the particular proximal gradient algorithm PANOC$^+$
	can be found in \cite{DeMarchiThemelis2022}.
An extension of the findings in \cite{KanzowMehlitz2022}, using a nonmonotone line search,
is given in \cite{DeMarchi2022}. In contrast to most existing papers on
proximal gradient methods, however, convergence of the entire
sequence is not addressed in \cite{DeMarchi2022,DeMarchiThemelis2022,KanzowMehlitz2022}.
Hence, no associated rate-of-convergence results could be given
(\cite{DeMarchi2022} presents some standard worst-case rate-of-convergence results
addressing the difference of two consecutive iterates along convergent
subsequences).
The aim of this paper is to fill this gap. 
More precisely, we show that the entire sequence 
generated by the proximal gradient method converges to 
a limit (with a suitable rate), provided that
this point 
	is an accumulation point of the generated sequence which
satisfies the Kurdyka--{\L}ojasiewicz property. 
The underlying convergence theory is still based on a merely local Lipschitz assumption on $\nabla f$,
neither its global Lipschitzness nor the (a priori)
boundedness of the iterates and stepsizes is presumed. 
To this end, we stress
that our analysis is not based on any kind of (global) descent lemma,  
which is in contrast to the contributions  
\cite{BauschkeBolteTeboulle2017,BolteSabachTeboulleVaisbourd2018} 
mentioned above.
	Let us emphasize that the mild assumptions used in 
	\cite{DeMarchi2022,DeMarchiThemelis2022,KanzowMehlitz2022}
	or the present paper do not guarantee the existence of accumulation points of the generated sequence.
	Some additional properties of the considered model problem are needed to ensure this.
	Since extensive numerical comparisons of (different types of) proximal gradient methods can already
	be found in several papers, see e.g.\
	\cite{GuWangHuoHuang2018,LiLin2015,WrightNowakFigueiredo2009},
	we abstain from the presentation of computational results here
	but focus on the compact justification of our theoretical findings.

The paper is organized as follows: 
In \cref{Sec:Prelims}, we formally introduce the model problem of interest and
provide some necessary notation as well as 
background material from generalized differentiation. 
The proximal gradient method together with the global convergence 
properties known from \cite{KanzowMehlitz2022} 
are stated in \cref{Sec:GenSpecGrad}.
The convergence and rate-of-convergence analysis is then given 
in \cref{Sec:Rate-of-Convergence}. We close with some final
remarks in \cref{Sec:Final}.

\section{Problem Setting and Preliminaries}\label{Sec:Prelims}

\subsection{Problem Setting}

Throughout the paper, we investigate the numerical treatment of
the composite optimization problem
\begin{equation}\label{Eq:P}\tag{P}
	\min_x \ \psi(x):=f(x) + \phi (x), \qquad x \in \mathbb X,
\end{equation}
where $ f\colon \mathbb{X} \to \mathbb{R} $ is
continuously differentiable, 
$ \phi\colon \mathbb{X} \to \overline{\mathbb{R}}:= \R\cup\{\infty\} $ 
is lower semicontinuous (possibly infinite-valued and nondifferentiable), 
and $ \mathbb{X} $ denotes a
Euclidean space, i.e., a real and finite-dimensional Hilbert space.
Since we do not want to deal with trivial situations, we assume that there exist
points in $\mathbb X$ where the value of $\phi$ is finite.
Let us underline that $\mathbb X$ is chosen to be Euclidean because this
allows to cover applications from matrix analysis like low-rank optimization
or matrix completion.

In order to minimize the function $\psi\colon\mathbb X\to\overline{\R}$ in \eqref{Eq:P},
we will exploit its composite structure which allows for gradient steps with respect 
to the continuously differentiable function $f$ on the one hand and so-called
proximal steps with respect to $\phi$ on the other hand, i.e., we rely on a
splitting approach. Throughout the last decades, experiments on numerous practically
relevant optimization problems have shown that splitting methods are superior
to the direct applications of standard methods from nonsmooth optimization to the
function $\psi$.

\subsection{Basic Notation}

Throughout the paper, the Euclidean space $\mathbb X$ 
will be equipped with the inner product 
$\langle\cdot,\cdot\rangle\colon \mathbb X\times\mathbb X\to\R$
and the associated norm $\norm{\cdot}$. Given a set 
$ A \subset \mathbb X $ and an element $ x \in \mathbb X $, 
we use $ A + x := x + A := \{ x \} + A :=
\{ x + a \mid a \in A \} $ for brevity.
Furthermore,
\[
	\dist(x,A):=\inf\{\norm{y-x}\,|\,y\in A\}
\]
denotes the distance of the point $x$ to the set $A$ with
$\dist(x,\emptyset):=\infty$.
For given $\varepsilon>0$, $B_\varepsilon(x):=\{y\in\mathbb X\,|\,\norm{y-x}\leq\varepsilon\}$
denotes the closed $\varepsilon$-ball around $x$.

The continuous linear operator
$f'(x)\colon\mathbb X\to\R$ denotes the derivative 
of the continuously differentiable function $f\colon\mathbb X\to\R$ 
at $x\in\mathbb X$, and we will
make use of $\nabla f(x):=f'(x)^*1$ where $f'(x)^*\colon\R\to\mathbb X$ 
is the adjoint of $f'(x)$.
This way, $\nabla f$ is a mapping from $\mathbb X$ to $\mathbb X$.

We further say that a sequence $ \{ x^k \} \subset \mathbb X $
converges \emph{Q-linearly} to $ x^*\in\mathbb X $ if there is a constant
$ c \in (0,1) $ such that the inequality
\[
   \| x^{k+1} - x^* \| \leq c \| x^k - x^* \|
\]
holds for all sufficiently large $ k \in \N $. Furthermore,
$\{x^k\}$ is said to converge \emph{R-linearly} to $ x^* $ if we
have
\[
   \limsup_{k \to \infty} \| x^k - x^* \|^{1/k} < 1.
\]
Note that this R-linear convergence holds if there 
exist constants $ \omega > 0 $ and $ \mu \in (0,1) $ such that
$ \| x^k - x^* \| \leq \omega \mu^k $ holds for all sufficiently
large $ k\in\N $, i.e., if the expression $ \| x^k - x^* \| $ is 
dominated by a Q-linearly convergent null sequence.

\subsection{Generalized Differentiation}
 
The following concepts are standard in variational analysis, and
we refer the interested reader to the monographs 
\cite{Mordukhovich2018,RockafellarWets2009} for more details.

Let us fix a merely lower semicontinuous function $\vartheta\colon\mathbb X\to\overline\R$
and pick $x\in\dom\vartheta$ where $\dom\vartheta:=\{x\in\mathbb X\,|\,\vartheta(x)<\infty\}$
denotes the domain of $\vartheta$. Then the set
\[
	\widehat\partial\vartheta(x)
	:=
	\left\{
		\eta\in\mathbb X\,\middle|\,
		\liminf\limits_{y\to x,\,y\neq x}
		\frac{\vartheta(y)-\vartheta(x)-\langle \eta,y-x\rangle}{\norm{y-x}}\geq 0
	\right\}
\] 
is called the {\em regular} (or {\em Fr\'{e}chet}) {\em subdifferential} of $\vartheta$ at $x$.
Furthermore, the set
\[
	\partial\vartheta(x)
	:=
	\left\{
		\eta\in\mathbb X\,\middle|\,
		\begin{aligned}
		&\exists\{x^k\},\{\eta^k\}\subset\mathbb X\colon\\
		&\qquad
		x^k\to x,\,\vartheta(x^k)\to\vartheta(x),\,\eta^k\to\eta,\,
		\eta^k\in\widehat{\partial}\vartheta(x^k)\,\forall k\in\N
		\end{aligned}
	\right\}
\]
is well known as the {\em limiting} (or {\em Mordukhovich}) 
{\em subdifferential} of $\vartheta$ at $x$.
Clearly, we always have $\widehat{\partial}\vartheta(x)\subset\partial\vartheta(x)$ 
by construction of these sets.
Whenever $\vartheta$ is a convex function, equality holds, 
and both subdifferentials coincide with the
subdifferential of convex analysis, i.e.,
\[
	\widehat{\partial}\vartheta(x)
	=
	\partial\vartheta(x)
	=
	\{
		\eta\in\mathbb X\,|\,
		\forall y\in\dom\vartheta\colon\,\vartheta(y)\geq\vartheta(x)+\langle\eta,y-x\rangle
	\}
\]
is valid in this situation.
By definition of the regular subdifferential, it is clear that
whenever $x^*\in\dom\vartheta$ is a local minimizer
of $\vartheta$, then $0\in\widehat\partial\vartheta(x^*)$ hold.
The latter fact is known as Fermat's rule, see \cite[Proposition~1.30(i)]{Mordukhovich2018}.
Thus, the inclusion $ 0 \in \partial \vartheta (x^*) $ is
a necessary optimality condition for $ x^* $ being a local
minimizer of $ \vartheta $ as well.
Note that, for $\vartheta $ being convex, this necessary optimality condition is also
sufficient for (global) minimality of $x^*$ for $\vartheta$.

Let us now apply this to the special case where $\vartheta:=\psi$ is the sum of the
continuously differentiable function $f$ and a merely lower semicontinuous function $\phi$,
as it happens to be the case when investigating \eqref{Eq:P}.
Whenever $x\in\dom\phi$ is fixed, the sum rule
\begin{equation}\label{eq:sum_rule}
		\partial (f+\phi)(x)
		=
		\nabla f(x)+ \partial \phi(x)
\end{equation}
holds due to the assumed continuous differentiability of $ f $, 
see \cite[Proposition~1.30(ii)]{Mordukhovich2018}.
Application of Fermat's rule therefore shows that the optimality
condition
\begin{equation*}
	0 \in \nabla f (x^*)+ \partial \phi(x^*)
\end{equation*}
holds at any local minimizer $ x^*\in\dom\phi $ of the composite optimization
problem \eqref{Eq:P}. Any point $ x^* \in \dom \phi $ satisfying
this necessary optimality condition will be called an 
\emph{M-stationary point} of \eqref{Eq:P}
due to the appearance of the limiting (or Mordukhovich) subdifferential.

We next introduce the famous Kurdyka--{\L}ojasiewicz property that 
was already mentioned in \cref{Sec:Intro} and which plays a 
central role in our subsequent convergence analysis. The 
version of this property stated below is a generalization of the classical 
Kurdyka--{\L}ojasiewicz inequality for nonsmooth functions as
introduced in
\cite{AttouchBolteRedontSoubeyran2010,BolteDaniilidisLewis2007,BolteDaniilidisLewisShiota2007}
and afterwards used in the local convergence analysis of several
nonsmooth optimization methods, see
\cite{attouch2009convergence,AttouchBolteSvaiter2013,BolteSabachTeboulle2014,BotCsetnek2016,BotCsetnekLaszlo2016,Ochs2018,OchsChenBroxPock2014}
for a couple of examples.

\begin{definition}\label{Def:KL-property}
Let $ g\colon \mathbb X \to \overline{\R} $ be 
lower semicontinuous. We say that $ g $ has the 
\emph{KL property}, where KL abbreviates \emph{Kurdyka--{\L}ojasiewicz}, 
at $ x^* \in \{ x \in \mathbb X \,|\, \partial g(x)\neq\emptyset \} $ if there exist a constant 
$ \eta > 0 $, a neighborhood $ U\subset\mathbb X $ of $ x^* $, and a continuous
concave function $ \chi\colon[0, \eta] \to [0,\infty) $ 
which is continuously differentiable on $(0,\eta)$ and satisfies $\chi(0)=0$ as well as
$\chi'(t)>0$ for all $t\in(0,\eta)$
such that the so-called \emph{KL inequality}
\[
   \chi ' \big( g(x) - g(x^*) \big) \dist \big( 0, \partial 
   g(x) \big) \geq 1
\]
holds for all $ x \in U \cap \big\{ x \in \mathbb X \,|\, 
g(x^*) < g(x) < g(x^*) + \eta \big\} $.
The function $ \chi $ from above is referred to as the \emph{desingularization function}.
\end{definition}

We note that there exist classes of functions where the 
KL property holds with the corresponding desingularization function
given by $ \chi (t) := c t^{\kappa} $ for 
$ \kappa \in (0,1] $ and some constant $ c > 0 $, where the
parameter $ \kappa $ is called the \emph{KL exponent},
see \cite{BolteDaniilidisLewisShiota2007,Kurdyka1998}.

\section{A Proximal Gradient Method and its Global Convergence Properties}\label{Sec:GenSpecGrad}

This section begins with a formal description of a 
proximal gradient method for the composite optimization
problem \eqref{Eq:P}, and then summarizes the associated
global convergence properties established in \cite{KanzowMehlitz2022}.
Note that our proximal gradient method uses a line
search which is important to get global convergence
properties without a global Lipschitz assumption. 
We start with a precise statement of the algorithm.

\begin{algorithm}[Proximal Gradient Method]\leavevmode
	\label{Alg:MonotoneProxGrad}
	\begin{algorithmic}[1]
		\REQUIRE $\tau > 1$, $0 < \gamma_{\min} \leq  \gamma_{\max} < \infty$, 
			$\delta \in (0,1)$, $x^0 \in \dom\phi $
		\STATE Set $k := 0$.
		\WHILE{A suitable termination criterion is violated at iteration $ k $}
		\STATE Choose $ \gamma_k^0 \in [ \gamma_{\min}, \gamma_{\max}] $.
		\STATE\label{step:subproblem_solve_MonotoneProxGrad} 
			For $ i = 0, 1, 2, \ldots $, compute a solution $ x^{k,i} $ of
      		\begin{equation}\label{Eq:Subki}
         		\min_x \ f (x^k) + \langle\nabla f(x^k), x - x^k \rangle + \frac{\gamma_{k,i}}{2} \| x - x^k \|^2 + \phi (x), 
         		\quad x \in \mathbb X
      		\end{equation}
      		with $ \gamma_{k,i} := \tau^i \gamma_k^0 $, until the acceptance criterion
      		\begin{equation}\label{Eq:StepCrit}
         		\psi (x^{k,i}) \leq 
         		\psi (x^k) - \delta \frac{\gamma_{k,i}}{2} \| x^{k,i} - x^k \|^2 
      		\end{equation}
      		holds.
		\STATE Denote by $ i_k := i $ the terminal value, and set $ \gamma_k := 
      			\gamma_{k,i_k} $ and $ x^{k+1} := x^{k,i_k} $.
      	\STATE Set $ k \leftarrow k + 1 $.
		\ENDWHILE
		\RETURN $x^k$
	\end{algorithmic}
\end{algorithm}

Our convergence analysis requires some technical assumptions
as well as a local Lipschitz condition on the gradient of 
the continuously differentiable function $f$.

\begin{assumption}\label{Ass:ProxGradMonotone}
\leavevmode
\begin{enumerate}
   \item \label{item:psi_bounded} The function $ \psi $ is bounded from below on $ \dom\phi $.
   \item \label{item:phi_bounded_affine} The function  $ \phi $ is bounded from below by an affine function.
   \item \label{item:local_Lipschitz} The function
   $ \nabla f\colon \mathbb X \to \mathbb X $ is locally Lipschitz continuous.
\end{enumerate}
\end{assumption}

Keeping in mind that our goal is to minimize the function $\psi$ in \eqref{Eq:P},
\cref{Ass:ProxGradMonotone}~\ref{item:psi_bounded} is reasonable.
Furthermore, \cref{Ass:ProxGradMonotone}~\ref{item:phi_bounded_affine}
is employed to guarantee existence of solutions for the appearing subproblems \eqref{Eq:Subki}.
To be precise, \Cref{Ass:ProxGradMonotone}~\ref{item:phi_bounded_affine} implies 
that the objective function of the subproblem \eqref{Eq:Subki} is, for fixed 
$ k, i \in\N$, coercive, and therefore always
attains a global minimizer $ x^{k,i} $ (which does not need to be unique). 
Finally, the local Lipschitz condition for $\nabla f$ from 
\Cref{Ass:ProxGradMonotone}~\ref{item:local_Lipschitz} will 
play a crucial role especially in \cref{Sec:Rate-of-Convergence} where
we consider situations where a sequence generated by \cref{Alg:MonotoneProxGrad}
converges as a whole and give associated rate-of-convergence results.

In the following, we recall the central global convergence
properties of \cref{Alg:MonotoneProxGrad} whose proofs can
be found in \cite[Section~3]{KanzowMehlitz2022}. Note that, throughout our analysis
of \cref{Alg:MonotoneProxGrad}, we implicitly assume that
this method generates an infinite sequence. For a discussion
of a practical termination criterion, we refer to 
\cite[Remark~3.1]{KanzowMehlitz2022} for more details.

First, we recall that the stepsize rule in \cref{step:subproblem_solve_MonotoneProxGrad} 
of \cref{Alg:MonotoneProxGrad} is always finite if the current iterate is not already
stationary. Hence, the overall method is well-defined.

\begin{lemma}\label{Lem:StepsizeFinite}
Consider a fixed iteration $ k\in\N $ of \cref{Alg:MonotoneProxGrad}, 
assume that $ x^k $ is not an M-stationary point of \eqref{Eq:P}, 
and suppose that \cref{Ass:ProxGradMonotone}~\ref{item:phi_bounded_affine} holds. 
Then the inner loop in \cref{step:subproblem_solve_MonotoneProxGrad} 
of \Cref{Alg:MonotoneProxGrad} is finite, i.e., we have $ \gamma_k = \gamma_{k,i_k} $
for some finite index $ i_k \in \{ 0, 1, 2, \ldots \} $.
\end{lemma}

The following result summarizes some of the properties of \Cref{Alg:MonotoneProxGrad} 
that will later be used in \cref{Sec:Rate-of-Convergence}. 

\begin{proposition}\label{Prop:xdiff}
Let \Cref{Ass:ProxGradMonotone}~\ref{item:psi_bounded} and \ref{item:phi_bounded_affine} hold,
and let $ \{ x^k \} $ be a sequence generated by \Cref{Alg:MonotoneProxGrad}.
Then the following statements hold:
\begin{enumerate}
	\item\label{item:diff_to_zero} $ \| x^{k+1} - x^k \| \to 0 $ as $ k \to \infty $,
	\item\label{item:stepsize_times_diff_to_zero} 
		for any convergent subsequence $ \{ x^k \}_K $, $ \gamma_k \| x^{k+1} - x^k \| \to_K 0 $ holds
		as $k\to_K\infty$,
	\item\label{item:stepsize_bounded}
		if, additionally, \cref{Ass:ProxGradMonotone}~\ref{item:local_Lipschitz} is valid,
		then for any convergent subsequence $\{ x^k \}_K$, $\{\gamma_k\}_K$ is bounded.
\end{enumerate}
\end{proposition}

Finally, we restate the main global convergence result for \cref{Alg:MonotoneProxGrad}, 
see again \cite[Section~3]{KanzowMehlitz2022} for the corresponding details.

\begin{theorem}\label{Thm:ConvProxGrad}
Let \Cref{Ass:ProxGradMonotone} be satisfied.
Then each accumulation point
of a sequence $ \{ x^k \} $ generated by \Cref{Alg:MonotoneProxGrad}
is an M-stationary point of \eqref{Eq:P}.
\end{theorem}

Note that \cite[Theorem~3.1]{KanzowMehlitz2022} shows that a result
like \Cref{Thm:ConvProxGrad} also holds without any Lipschitz
condition regarding $ \nabla f $, but it then requires
a slightly stronger condition for the nonsmooth function $ \phi $,
namely the continuity of $ \phi $ on its domain 
(this condition holds, e.g., if $ \phi $ is the indicator function of a constraint set). 
Our analysis in \cref{Sec:Rate-of-Convergence},
however, requires the local Lipschitz condition for the 
gradient $ \nabla f $, so we decided to treat it as a standing assumption.

We close this section by mentioning two classes of examples
where the standard global Lipschitz assumption on the
gradient of $  f $ is typically violated, whereas a
local Lipschitz condition is often satisfied.

\begin{example}\label{ex:alm} (Augmented Lagrangian Methods) \\
Consider the constrained optimization problem 
\[
   \min\limits_x \ f(x) + \phi (x) \quad \text{s.t.} \quad c(x) \in C,
\]
where $ f\colon\mathbb X\to\R$ and $\phi\colon\mathbb X\to\overline\R$ 
are as in \eqref{Eq:P}. 
In addition, we have some constraints defined by a continuously differentiable function 
$ c\colon\mathbb X\to \mathbb Y$, where $\mathbb Y$ is another Euclidean space, and a 
nonempty, closed, and convex set $ C\subset\mathbb Y $. 

Given a current iterate $ x^k\in\mathbb X $ and a corresponding Lagrange multiplier
estimate $ \lambda^k\in\mathbb Y $, augmented Lagrangian techniques then
compute the next iterate $ x^{k+1} $ by solving (approximately)
the subproblem
\[
   \min_x 
   f(x) + \phi (x) 
   + 
   \frac{\rho_k}{2} 
   \dist^2 
  	 \left( c(x) + \frac{\lambda^k}{\rho_k},C \right),
   \qquad  
   x \in \mathbb X
\]
for some penalty parameter $ \rho_k > 0 $. Since the squared
distance function $y\mapsto\dist^2(y,C)$ is continuously differentiable by convexity of $C$,
see \cite[Corollary~12.31]{BauschkeCombettes2017}, this subproblem has exactly the
structure of the composite optimization problem \eqref{Eq:P}
and can therefore, in principle, be solved by a proximal
gradient method, 
see \cite{ChenGuoLuYe2017,GuoDeng2021,JiaKanzowMehlitzWachsmuth2021,DeMarchiJiaKanzowMehlitz2022} 
for suitable realizations of this approach. 

Assuming that the gradient of the smooth part of this objective function
(with respect to the variable $x$)
is globally Lipschitz continuous, however, is pretty strong
is this setting and, basically, requires the constraint function
$ c $ to be linear and the set $ C $ to be polyhedral, whereas
local Lipschitzness of this gradient holds under mild
conditions on the smoothness of $f$ and $c$.
\end{example}

The following example makes use of conjugate functions, 
see \cite[Definition~13.1]{BauschkeCombettes2017}.
Since, within this paper, they only occur in this particular
application, we refrain from stating their precise definitions and
properties, and refer the interested reader to the excellent
monographs \cite{BauschkeCombettes2017,Beck2017,RockafellarWets2009} for more details.

\begin{example}\label{ex:dual_prox_grad} (Dual Proximal Gradient Methods) \\
Consider the (primal) optimization problem 
\begin{equation}\label{Eq:Primal}
   \min\limits_x \ g(x) + h(Ax), \qquad x \in \mathbb{X}
\end{equation}
where both functions $g\colon\mathbb X\to\overline{\R}$ and
$h\colon\mathbb Y\to\overline{\R}$ are lower semicontinuous
and convex while possessing nonempty domains, 
and $ A\colon\mathbb X\to\mathbb Y$ 
is a linear operator. 
Above, $\mathbb Y$ is another Euclidean space.
Note that none of the functions $g$ or $h$ is assumed to be
(continuously) differentiable.

 The (Fenchel) dual problem of \eqref{Eq:Primal} is given by
\begin{equation}\label{Eq:Dual}
	\min\limits_y \ g^* (A^* y) + h^* (-y),\qquad y\in\mathbb Y
\end{equation} 
with the two conjugate functions $ g^*\colon\mathbb X\to\overline{\R}$ 
and $h^*\colon\mathbb Y\to\overline{\R} $ being lower semicontinuous and convex,
and $A^*\colon\mathbb Y\to\mathbb X$ being the adjoint of $A$.
Under suitable assumptions, the pair \eqref{Eq:Primal}, \eqref{Eq:Dual}
enjoys strong duality, i.e., the optimal objective function values of these
problems coincide, see \cite{Rockafellar1970}, which motivates to solve
\eqref{Eq:Dual} instead of \eqref{Eq:Primal} in some applications where the
conjugate functions are explicitly available.

Assuming, in addition, 
that $ g $ is uniformly convex, it is known that $ g^* $ is
real-valued everywhere and continuously differentiable with
a globally Lipschitz continuous gradient, 
see \cite[Proposition~12.60]{RockafellarWets2009}.
Consequently, as promoted in \cite{BeckTeboulle2014}, 
a standard proximal gradient algorithm can be applied to the dual problem \eqref{Eq:Dual}. 
On the other hand, if $ g $ is only strictly convex, then the domain of $ g^* $ is,
in general, no longer the entire space, but $ g^* $ can still
be shown to be continuously differentiable on the interior of
its domain. Its gradient, however, is no longer guaranteed to
be globally Lipschitz continuous on the domain.
\end{example}

\section{Convergence Analysis in the Presence of the KL Property}\label{Sec:Rate-of-Convergence}

The aim of this section is to show convergence of the entire 
sequence $ \{ x^k \} $ generated by \Cref{Alg:MonotoneProxGrad}
provided that there exists an accumulation point $ x^* $ which,
in addition, satisfies the KL property, and to present associated
rate-of-convergence results. The proofs of these results are based
on a local Lipschitz assumption on $\nabla f$ only, without
the a priori assumption that the whole sequence $ \{ x^k \} $ is 
bounded. Based on some recent contributions in the area of
proximal gradient and related first-order methods, it seems
reasonable to expect such a result to hold. For example,
\cite{BolteSabachTeboulle2014,Ochs2018} consider a whole class of first-order methods 
and investigate their (essentially local) convergence showing, in particular,
that the entire sequence $ \{ x^k \} $ generated by their methods stays 
within a certain neighborhood of a solution
provided that the KL property holds at this solution. 

Their approach is not directly applicable to our situation since,
on the one hand, we do not use the a priori assumption that
our iterates are bounded, and, on the other hand, 
because the adaption of the methods considered in 
\cite{BolteSabachTeboulle2014,Ochs2018} to the 
proximal gradient setting would result in an algorithm with a 
constant stepsize. However, having an accumulation point of \Cref{Alg:MonotoneProxGrad}
satisfying the KL property, we know from the local Lipschitz 
assumption on $\nabla f$ that a respective global Lipschitz condition holds in a 
suitable neighborhood of this point, which then can be used to
verify that the stepsizes computed by \Cref{Alg:MonotoneProxGrad}
remain bounded.
This -- more or less heuristic -- idea fortifies us to believe 
that one can also get convergence and rate-of-convergence
results under the KL property in the presence of 
\cref{Ass:ProxGradMonotone}~\ref{item:local_Lipschitz}. 
The following analysis is a careful mathematical realization of this somewhat vague idea.

We begin with a result which shows that, locally around an 
accumulation point of the sequence $\{x^k\}$, 
the associated stepsizes $ \gamma_k $ remain bounded.
This observation and its proof are related to 
\cite[Corollary~3.1]{KanzowMehlitz2022}.
Note that this statement is essentially different from the boundedness of stepsizes
along convergent subsequences of iterates which is inherent in the presence of
\cref{Ass:ProxGradMonotone}, see \cref{Prop:xdiff}~\ref{item:stepsize_bounded}.

\begin{lemma}\label{Lem:gamma-bound}
Let \Cref{Ass:ProxGradMonotone} hold, let $ \{ x^k \} $
be any sequence generated by \Cref{Alg:MonotoneProxGrad},
and let $ x^* $ be an accumulation point of this sequence.
Then, for any $ \rho > 0 $, there is
a constant $ \bar \gamma_{\rho} > 0 $ (usually depending 
on $ \rho $) such that $ \gamma_k \leq \bar \gamma_{\rho} $
holds for all $ k \in \N $ such that
$ x^k \in B_{\rho} (x^*) $.
\end{lemma}

\begin{proof}
First, recall from \Cref{Lem:StepsizeFinite} that the 
stepsize $ \gamma_k $ is well-defined for each $k\in\N$. 
Let $ \rho > 0 $ be fixed, and recall that the assumed
local Lipschitz continuity of $ \nabla f $ implies that 
this gradient mapping is (globally) Lipschitz continuous
on the compact set $ B_{2 \rho} (x^*) $ (note that we
took $ 2 \rho $ as the radius of this ball here). 
Let us denote the corresponding Lipschitz constant 
by $ L_{2\rho} $. Since $ x^* $ is an accumulation point of the sequence $ \{ x^k \} $, 
there are infinitely many iterates of this sequence belonging to $ B_{\rho} (x^*) $. 

Now, assume, by contradiction, that there is a subsequence
$ \{ \gamma_k \}_K $ with $ x^k \in B_{\rho} (x^*) $ for all
$ k \in K $ such that $ \{ \gamma_k \}_K $ is unbounded. 
Without loss of generality, we may assume that $ \gamma_k \to_K \infty $,
that the subsequence of iterates $ \{ x^k \}_K $ converges 
to some point $ \bar x $ (not necessarily equal to $ x^* $),
and that, for each $ k \in K $,
the acceptance criterion \eqref{Eq:StepCrit} is violated 
in the first iteration of the inner loop. 
Then, for the trial stepsize 
$ \hat \gamma_k := \gamma_k / \tau = \tau^{i_k-1}\gamma_k^0$, 
we also have $ \hat \gamma_k \to_K \infty $, 
whereas the corresponding trial vector $ \hat x^k := x^{k, i_k-1} $
does not satisfy the acceptance criterion from \eqref{Eq:StepCrit},
i.e., we have
\begin{equation}\label{Eq:3-12}
	\psi (\hat x^k) > \psi (x^k) - \delta \frac{\hat \gamma_k}{2}
	\| \hat x^k - x^k \|^2 \quad \forall k \in K.
\end{equation}
On the other hand, since $ \hat x^k $ solves the corresponding
subproblem \eqref{Eq:Subki} with $ \hat \gamma_k $ in place of $\gamma_{k,i}$, we have
\begin{equation}\label{Eq:3-13}
	\langle \nabla f(x^k), \hat x^k - x^k \rangle + 
	\frac{\hat \gamma_k}{2} \| \hat x^k - x^k \|^2 +
	\phi (\hat x^k) - \phi (x^k) \leq 0.
\end{equation}
We claim that this, in particular, implies that $ \hat x^k \to_K \bar x $. 
In fact, using \eqref{Eq:3-13}, the Cauchy-Schwarz inequality,
and the fact that $ \{ \psi (x^k) \} $ is monotonically decreasing
by construction of \cref{Alg:MonotoneProxGrad}, we obtain
\begin{align*}
	\frac{\hat \gamma_k}{2} \| \hat x^k - x^k \|^2 & \leq 
	\| \nabla f(x^k) \| \|  \hat x^k - x^k \| + \phi (x^k) -
	\phi (\hat x^k ) \\ 
	& = \| \nabla f(x^k) \| \|  \hat x^k - x^k \| + \psi (x^k)
	- f(x^k) - \phi (\hat x^k) \\
	& \leq \| \nabla f(x^k) \| \|  \hat x^k - x^k \| +
	\psi (x^0) - f(x^k) - \phi (\hat x^k).
\end{align*}
Since $ f $ is continuously differentiable and $ - \phi $ is 
bounded from above by an affine function in view of 
\Cref{Ass:ProxGradMonotone}~\ref{item:phi_bounded_affine}, the above estimate
implies $ \| \hat x^k - x^k \| \to_K 0 $. In fact, if 
$ \{ \| \hat x^k - x^k \| \}_K $ would be unbounded, then the 
left-hand side would grow more rapidly than the right-hand 
side, and if $ \{ \| \hat x^k - x^k \| \}_K $ would be bounded,
but staying away, at least on a subsequence, from zero by a 
positive number, the right-hand side would be bounded, whereas
the left-hand side would be unbounded on the corresponding
subsequence. Consequently, we have $ \| \hat x^k - x^k \| \to_K 0 $, 
and since $ x^k \to_K \bar x $, this implies
$\hat x^k \to_K \bar x $. In particular, since 
$ \bar x \in B_{\rho} (x^*) $, this implies that, for all
sufficiently large $ k \in K $, we have both $ x^k \in 
B_{2\rho} (x^*) $ and $ \hat x^k \in B_{2 \rho} (x^*) $.

Let us fix some $k\in K$.
Using the mean-value theorem yields the existence of a point $ \xi^k $ 
on the line segment connecting $ x^k $ with $ \hat x^k $ such that 
\begin{align*}
	\psi (\hat x^k) - \psi (x^k) 
	&= 
	f(\hat x^k) + \phi (\hat x^k)-f(x^k) - \phi (x^k) 
	\\
	&= 
	\langle \nabla f(\xi^k) , 	\hat x^k - x^k \rangle + \phi (\hat x^k) - \phi (x^k).
\end{align*}
Substituting the resulting expression for $ \phi (\hat x^k) - \phi (x^k) $
into \eqref{Eq:3-13}, we see that
\begin{equation}\label{Eq:3-14a}
	\langle \nabla f(x^k) - \nabla f(\xi^k) ,
	\hat x^k - x^k \rangle + \frac{\hat \gamma_k}{2}
	\| \hat x^k - x^k \|^2 + \psi (\hat x^k) - \psi (x^k) \leq 0. 
\end{equation}
Exploiting \eqref{Eq:3-12}, we therefore obtain 
\begin{align*}
	\frac{\hat \gamma_k}{2} \| \hat x^k - x^k \|^2 
	& \leq - \langle \nabla f(x^k) - \nabla f(\xi^k) ,
	\hat x^k - x^k \rangle + \psi (x^k) - \psi (\hat x^k) \\
	& \leq \| \nabla f(x^k) - \nabla f(\xi^k ) \| \| \hat x^k -
	x^k \| + \delta \frac{\hat \gamma_k}{2} \| \hat x^k - x^k
	\|^2 
\end{align*}
which can be rewritten as
\[
	( 1 - \delta ) \frac{\hat \gamma_k}{2} \| \hat x^k - x^k \| 
	\leq \| \nabla f(x^k) - \nabla f (\xi^k) \|.
\]
Since $ \xi^k $ in an element from the line connecting
$ x^k $ and $ \hat x^k $, it follows that $ \xi^k \in 
B_{2 \rho} (x^*) $ for all $ k \in K $ sufficiently large.
Hence, the Lipschitz continuity of $ \nabla f $ on this
ball, we find
\begin{equation*}
	( 1 - \delta ) \frac{\hat \gamma_k}{2} \| \hat x^k - x^k \| 
	\leq L_{2\rho} \| x^k - \xi^k \| \leq L_{2\rho} 
	\| x^k - \hat x^k \| 
\end{equation*}
for all sufficiently large $k\in K$.
Since $ \hat x^k \neq x^k $ in view of \eqref{Eq:3-12}, 
this implies that $ \{ \hat \gamma_k \}_K $ is bounded which,
in turn, yields the boundedness of the subsequence
$ \{ \gamma_k \}_K $, contradicting our assumption.
This completes the proof.
\end{proof}

We next show  that the entire 
sequence $ \{ \psi (x^k )\} $
converges to $ \psi (x^*) $, where $ x^* $ is an 
arbitrary accumulation point of a sequence $ \{ x^k \} $
generated by \Cref{Alg:MonotoneProxGrad}. Note that this result
is not completely obvious since $ \psi $ is only lower semicontinuous but
not continuous in general.
Indeed, this property results from the construction of the iterates $ x^{k+1} $ of
\Cref{Alg:MonotoneProxGrad}.

\begin{lemma}\label{Lem:psi-conv}
Let \Cref{Ass:ProxGradMonotone} be satisfied, and
let $ x^* $ be an accumulation point of a sequence $ \{ x^k \} $ generated by \Cref{Alg:MonotoneProxGrad}. 
Then the entire sequence $ \{ \psi (x^k) \} $ converges to $ \psi (x^*) $.
\end{lemma}

\begin{proof}
Let $ \{ x^k \}_K $ be a subsequence converging to $ x^* $.
By means of \Cref{Prop:xdiff}~\ref{item:diff_to_zero},
we also have $x^{k+1}\to_K x^*$.
Since $ \psi $ is lower semicontinuous, we then obtain 
\begin{equation}\label{eq:psi-lower-bound}
	\psi (x^*) \leq \liminf_{k \to_K\infty} \psi (x^{k+1}).
\end{equation}
On the other hand, by construction, the entire sequence 
$ \{ \psi (x^k) \} $ is monotonically decreasing. Since 
it is also bounded from below by $ \psi (x^*) $ as a
consequence of \eqref{eq:psi-lower-bound}, it follows that
the whole sequence $ \{ \psi (x^k) \} $ converges. It 
remains to show that its limit is equal to (the lower
bound) $ \psi (x^*) $.
	
To this end, we first note that
$ x^{k+1} $ solves the subproblem 
\eqref{Eq:Subki} with stepsize $ \gamma_k $
in place of $\gamma_{k,i}$.
Hence, we have 
\begin{align*}
	\langle \nabla f (x^k), x^{k+1} - x^k \rangle 
	&+
	\frac{\gamma_k}{2} \| x^{k+1} - x^k \|^2 + \phi (x^{k+1})
	\\
	&\leq 
	\langle \nabla f (x^k), x^* - x^k \rangle 
	+
	\frac{\gamma_k}{2} \| x^* - x^k \|^2 + \phi (x^*)
\end{align*}
for each $ k \in \N $. Taking the upper limit as $ k \to_K\infty$, 
and using the continuity of $ \nabla f $ as well as \Cref{Prop:xdiff},
we obtain 
\[
   \limsup_{k \to_K \infty} \phi (x^{k+1}) \leq \phi (x^*).
\]
Combining this with \eqref{eq:psi-lower-bound} and 
using the continuity of $f$ yields $\psi(x^{k+1})\to_K\psi(x^*)$.
Since $\{\psi(x^k)\}$ converges, the assertion follows. 
\end{proof}

All results stated so far are independent of the KL property.
The remaining part of our analysis, however, is heavily based on 
the assumption that our objective function $ \psi $
satisfies the KL property at a given accumulation point $ x^* $
of a sequence $\{x^k\}$ generated by \cref{Alg:MonotoneProxGrad}.
In particular, let $ \eta > 0 $ be the corresponding 
constant from the definition of the associated desingularization function $\chi$.
Furthermore, we will assume that \cref{Ass:ProxGradMonotone} is valid.
In view of \Cref{Prop:xdiff}, we can find a 
sufficiently large index $ \hat{k} \in \N $ such that
\begin{equation}\label{Eq:eta}
   \sup_{k \geq \hat{k}} \| x^{k+1} - x^k \| \leq \eta .
\end{equation}
We then define
\begin{equation}\label{Eq:rho}
   \rho := \eta + \frac{1}{2}
\end{equation}
as well as the compact set
\begin{equation}\label{Eq:C_rho}
	C_\rho := B_\rho(x^*)\cap\mathcal L_\psi(x^0),
\end{equation}
where $\mathcal L_\psi(x^0):=\{x\in\mathbb X\,|\,\psi(x)\leq\psi(x^0)\}$
is the sublevel set of $\psi$ with respect to $x^0$, the starting
point exploited in \cref{Alg:MonotoneProxGrad}.
By monotonicity of $\{\psi(x^k)\}$, we have $\{x^k\}\subset\mathcal L_\psi(x^0)$.
Finally, throughout the section, 
let $L_\rho>0$ be a (global) Lipschitz constant of $ \nabla f $ on $C_\rho$ from \eqref{Eq:C_rho}.
Finally, in view of \Cref{Lem:gamma-bound}, we have
\begin{equation}\label{Eq:gamma-rho}
   \gamma_k \leq \bar \gamma_{\rho} 
   \quad\forall x^k \in C_\rho
\end{equation}
with some suitable upper bound $ \bar \gamma_{\rho} > 0 $ 
(depending on our choice of $ \rho $ from \eqref{Eq:rho}).
Using this notation, we can formulate the following result.

\begin{lemma}\label{Lem:alpha-small}
Let \Cref{Ass:ProxGradMonotone} hold, and let $ \{ x^k \} $
be any sequence generated by \Cref{Alg:MonotoneProxGrad}.
Suppose that $ \{ x^k \}_K $ is a subsequence converging to
some limit point $ x^* $, and that $ \psi $ has the KL property
at $ x^* $ with desingularization function $ \chi $.
Then there is a sufficiently large constant $ k_0 \in K $
such that the corresponding constant 
\begin{equation}\label{Eq:alpha}
   \alpha := 
   \| x^{k_0} - x^* \| 
   +
   \sqrt{\frac{8 \big( \psi(x^{k_0}) - \psi (x^*) \big)}{\delta \gamma_{\min}}} 
   + 
   \frac{2 \big( \bar \gamma_{\rho} + 
   L_{\rho} \big)}{\delta \gamma_{\min}} \chi \big( 
   \psi (x^{k_0}) - \psi (x^*) \big)
\end{equation}
satisfies $ \alpha < \frac{1}{2} $, where $ \rho>0$ and $\bar\gamma_\rho>0$ are
the constants defined in \eqref{Eq:rho} and \eqref{Eq:gamma-rho}, respectively,
while $L_\rho>0$ is a Lipschitz constant of $\nabla f$ on $C_\rho$ from
\eqref{Eq:C_rho}, and
$ \delta>0 $ as well as $ \gamma_{\min}>0 $ are the parameters from
\Cref{Alg:MonotoneProxGrad}.
\end{lemma}

\begin{proof}
The statement follows from the fact that each summand on the 
right-hand side of \eqref{Eq:alpha} can be made arbitrarily
small. This is clear for the first one since the 
subsequence $ \{ x^k \}_K $ converges to $ x^* $. 
This is also true for the second summand as a consequence of 
\Cref{Lem:psi-conv}. Finally, the third one can be
made arbitrarily small since we have $ \psi (x^k) \to 
\psi (x^*)  $ by \Cref{Lem:psi-conv}, taking into account
that the desingularization function $ \chi $ is 
continuous at the origin. Hence, the statement follows
by taking an index $ k_0 \in K $ sufficiently large.
\end{proof}

We next state another technical result.

\begin{lemma}\label{Lem:DistanceSubgrad}
Let \Cref{Ass:ProxGradMonotone} hold, and let $ \{ x^k \} $
be any sequence generated by \Cref{Alg:MonotoneProxGrad}.
Suppose that $ \{ x^k \}_K $ is a subsequence converging to
some limit point $ x^* $, and that $ \psi $ has the KL property
at $ x^* $ with desingularization function $ \chi $.
Then
\begin{equation*}
	\dist\big( 0, \partial \psi (x^{k+1}) \big) \leq 
	\big( \bar \gamma_{\rho} + L_{\rho} \big) \| x^{k+1} - 
	x^k \|
\end{equation*}
holds for all sufficiently large 
$ k \in \N $ such that $ x^k \in B_{\alpha} (x^*) $,
where $ \alpha<\tfrac12 $ denotes the constant from \eqref{Eq:alpha},
$ \bar \gamma_{\rho}>0$ is the constant from \eqref{Eq:gamma-rho},
and $L_{\rho}>0$ is the Lipschitz constant of $\nabla f$ on $C_\rho$
from \eqref{Eq:C_rho}.
\end{lemma}

\begin{proof}
For any $ k \in \N $, since $ x^{k+1} $ is a solution of
\eqref{Eq:Subki}, we obtain
\begin{equation*}
	0 \in \nabla f (x^k) + \gamma_k (x^{k+1}-x^k) + 
	\partial \phi (x^{k+1})
\end{equation*}
from the corresponding M-stationary condition. This implies
\begin{equation}\label{Eq:MstatImplies}
	\gamma_k ( x^k - x^{k+1} ) + \nabla f (x^{k+1}) -
	\nabla f (x^k) \in \nabla f (x^{k+1}) +
	\partial \phi (x^{k+1}) = \partial \psi (x^{k+1})
\end{equation}
for all $ k \in \N $, where we used the sum rule \eqref{eq:sum_rule}
for the limiting subdifferential. 

Now, take an arbitrary index $ k \in \N $ sufficiently large
such that $ x^k \in B_{\alpha} (x^*) $ and $ k \geq \hat{k} $,
where $ \hat{k} $ is the index from \eqref{Eq:eta}. 
In view of \eqref{Eq:rho} and
\Cref{Lem:alpha-small}, we have $ \alpha \leq \rho $.
Therefore, \Cref{Lem:gamma-bound} shows that 
\begin{equation}\label{Eq:Tog-gammak}
	\gamma_k \leq \bar \gamma_{\rho} . 
\end{equation}
Moreover, using \eqref{Eq:eta}, \eqref{Eq:rho}, and \Cref{Lem:alpha-small},
we get
\begin{equation*}
	\| x^{k+1} - x^* \| \leq \| x^{k+1} - x^k \| + \| x^k - x^* \|
	\leq \eta + \alpha \leq \rho.
\end{equation*}
Hence, $ x^k, x^{k+1} \in C_{\rho} $ holds
with the compact set $ C_{\rho} $ from \eqref{Eq:C_rho}.
Therefore, we have
\begin{equation*}
	\big\| \nabla f (x^{k+1}) - \nabla f (x^k ) \big\| 
	\leq L_{\rho} \| x^{k+1} - x^k \|
\end{equation*}
by definition of $L_\rho$.
Together with \eqref{Eq:MstatImplies} and \eqref{Eq:Tog-gammak}, 
we thus obtain
\begin{align*}
	\dist\big( 0, \partial \psi (x^{k+1}) \big) &
	\leq \big\| \gamma_k ( x^k - x^{k+1} ) + \nabla f (x^{k+1}) -
	\nabla f (x^k)	\big\| \\
	& \leq \gamma_k \| x^{k+1} - x^k \| + L_{\rho} \| x^{k+1} - x^k \| \\
	& \leq ( \bar \gamma_{\rho} + L_{\rho} \big) \| x^{k+1} - x^k \|
\end{align*}
for all $ k \in \N$ satisfying $k \geq \hat{k}$ and $ x^k \in B_{\alpha} (x^*) $.
\end{proof}

The following result shows that the entire sequence 
$ \{ x^k \} $, generated by \cref{Alg:MonotoneProxGrad}, 
already converges to one of its accumulation points $ x^* $ 
provided that the objective function $ \psi $
satisfies the KL property at this point.
The proof combines our previous results with a technique
used in \cite{BolteSabachTeboulle2014}.

\begin{theorem}\label{Thm:GlobConv}
Let \Cref{Ass:ProxGradMonotone} hold, and let $ \{ x^k \} $
be any sequence generated by \Cref{Alg:MonotoneProxGrad}.
Suppose that $ \{ x^k \}_K $ is a subsequence converging to
some limit point $ x^* $, and that $ \psi $ has the KL property
at $ x^* $. Then the entire sequence $ \{ x^k \} $ converges
to $ x^* $.
\end{theorem}

\begin{proof}
In view of \Cref{Lem:psi-conv}, we know that the whole
 sequence $ \{ \psi (x^k )\} $ is 
monotonically decreasing and converging to $ \psi (x^*) $.
This implies that
$ \psi (x^k) \geq \psi (x^*) $ holds for all $ k \in \N $.

Now, suppose we have $ \psi (x^k ) = \psi (x^*) $ for some
index $ k \in \N $. Then, by monotonicity, we also get
$ \psi (x^{k+1}) = \psi (x^*) $. Consequently, we obtain 
from \eqref{Eq:StepCrit} that
\begin{equation*}
	0 \leq \frac{\delta \gamma_{\min}}{2} \| x^{k+1} - x^k \|^2 
	\leq \psi (x^k) - \psi (x^{k+1}) = 0
\end{equation*}
and, thus, $ x^{k+1} = x^k $. Since, by assumption, 
the subsequence $ \{ x^k \}_K $ converges to $ x^* $, this
implies that $ x^k = x^* $ for all $ k \in \N $ sufficiently
large. In particular, we have convergence of the entire 
(eventually constant) sequence
$ \{ x^k \} $ to $ x^* $ in this situation.

For the remainder of this proof, we can therefore assume that 
$ \psi(x^k) > \psi (x^*) $ holds for all $ k \in \N $.
We then let $ \alpha \in (0, 1/2) $ be the constant from 
\eqref{Eq:alpha}, and $ k_0 \in K $ be the corresponding
iteration index which is used in the definition of $ \alpha $,
see \cref{Lem:alpha-small}. We then have
$ 0 < \psi (x^k) - \psi (x^*) \leq \psi (x^{k_0}) - \psi (x^*) $
for all $k \geq k_0 $. Without loss of generality, we may
also assume that $ k_0 \geq \hat{k} $ (the latter being the
index defined by \eqref{Eq:eta}) and that $ k_0 $ is 
sufficiently large to satisfy
\begin{equation}\label{Eq:k0-large}
	\psi (x^{k_0}) < \psi (x^*) + \eta.
\end{equation}
Let $\chi\colon[0,\eta]\to[0,\infty)$ be the desingularization function
which comes along with the validity of the KL property at $x^*$.
Due to $ \chi (0) = 0 $ and $ \chi '(t) > 0 $ for all $t\in(0,\eta)$, we obtain
\begin{equation}\label{Eq:chi-get}
	\chi \big( \psi (x^k ) - \psi (x^*)\big) \geq 0 \quad \forall k \geq k_0.
\end{equation}

We now claim that the following two statements hold for all 
$ k \geq k_0 $:
\begin{enumerate}
	\item \label{Item:Ind-1}
	   $ x^k \in B_{\alpha} (x^*) $, 
	\item \label{Item:Ind-2}
	  $ \| x^{k_0} - x^* \| + \sum_{i=k_0}^k \| x^{i+1} - x^i \| \leq \alpha $, which is equivalent to
	  \begin{equation}\label{Eq:Ind-2}
	  	 \sum_{i=k_0}^k \| x^{i+1} - x^i \| \leq 
	  	 \sqrt{\frac{8 \big( \psi(x^{k_0}) - \psi (x^*) \big)}{\delta \gamma_{\min}}} + \frac{2 \big( \bar \gamma_{\rho} + 
	  	 L_{\rho} \big)}{\delta \gamma_{\min}} \chi \big( 
	  	 \psi (x^{k_0}) - \psi (x^*) \big).
	  \end{equation}
\end{enumerate}
We verify these two statements jointly by induction. 
For $ k = k_0 $, statement \ref{Item:Ind-1} holds simply by the 
definition of $ \alpha $ in \eqref{Eq:alpha}. Furthermore, 
the acceptance criterion \eqref{Eq:StepCrit} together with the 
monotonicity of $ \{ \psi (x^k) \} $ implies
\begin{equation}\label{Eq:induktiv}
   \| x^{k_0+1} - x^{k_0} \| \leq \sqrt{\frac{2 \big( 
   \psi (x^{k_0}) - \psi (x^{k_0+1})	\big)}{\delta 
   \gamma_{\min}}} \leq \sqrt{\frac{2 \big( 
   	\psi (x^{k_0}) - \psi (x^*) \big)}{\delta 
   	\gamma_{\min}}}.
\end{equation}
In particular, this shows that \eqref{Eq:Ind-2} holds for 
$ k = k_0 $. Suppose that both statements hold for some $ k \geq k_0 $. Using the triangle inequality, the induction hypothesis, 
and the definition of $ \alpha $, we obtain
\[
   \| x^{k+1} - x^* \| \leq \sum_{i=k_0}^k \| x^{i+1} - x^i \| 
   + \| x^{k_0} - x^* \| \leq \alpha ,
\]
i.e., statement \ref{Item:Ind-1} holds for $ k+1 $ in place of $k$. The verification of 
the induction step for \ref{Item:Ind-2} is more involved.

To this end, first note that \eqref{Eq:k0-large} implies
\begin{equation}\label{Eq:3-14}
	\psi (x^*) < \psi (x^i) < \psi (x^*) + \eta \quad \forall
	i \geq k_0.
\end{equation}
Since $ \psi $ has the KL property at $ x^* $, we have
\begin{equation}\label{Eq:3-15}
	\chi ' \big( \psi (x^i) - \psi (x^*) \big) \dist
	\big( 0, \partial \psi (x^i) \big) \geq 1 \quad 
	\forall i \geq k_0.
\end{equation}
Since $ x^i \in B_{\alpha} (x^*) $ for all $ i \in\{ k_0, k_0 +1, 
\ldots, k \}$ by our induction hypothesis, we can apply 
\Cref{Lem:DistanceSubgrad} and obtain (after a simple index shift)
\begin{equation*}
	\text{dist} \big( 0, \partial \psi (x^i) \big) \leq 
	\big( \bar \gamma_{\rho} + L_{\rho} \big) \| x^i - x^{i-1} \|
	\quad \forall i \in \{k_0 +1, k_0 + 2, \ldots, k+1 \}.
\end{equation*}
In view of \eqref{Eq:3-15}, we therefore obtain
\begin{equation}\label{Eq:3-16}
	\chi ' \big( \psi (x^i) - \psi (x^*) \big) \geq
	\frac{1}{\big( \bar \gamma_{\rho} + L_{\rho} \big) 
	\| x^i - x^{i-1} \|}\quad \forall 
	i \in \{ k_0 +1, k_0 + 2, \ldots, k+1 \}.
\end{equation}
To simplify some of the subsequent formulas, we follow
\cite{BolteSabachTeboulle2014} and introduce
the short-hand notation
\[
   \Delta_{i,j} := \chi \big( \psi (x^i) - \psi (x^*) \big)
   - \chi \big( \psi (x^j) - \psi (x^*) \big)
\]
for $ i, j \in \N $. The assumed concavity of $ \chi $ then 
implies
\begin{equation}\label{Eq:3-17}
	\Delta_{i,i+1} \geq \chi ' \big( \psi (x^i) - \psi (x^*) \big)
	\big( \psi (x^i) - \psi (x^{i+1}) \big).
\end{equation}
Using \eqref{Eq:3-16}, \eqref{Eq:3-17}, and the acceptance 
criterion \eqref{Eq:StepCrit}, we therefore get
\begin{align*}
	\Delta_{i,i+1} & \geq 
	\chi ' \big( \psi (x^i) - \psi (x^*) \big)
	\big( \psi (x^i) - \psi (x^{i+1}) \big) \\
	& \geq \frac{\psi (x^{i}) - \psi (x^{i+1})}{ ( \bar \gamma_{\rho}
	+ L_{\rho}) \| x^i - x^{i-1} \|} 
	 \geq \frac{\delta \gamma_{\min}}{2 (\bar \gamma_{\rho}
	+ L_{\rho})} \frac{\| x^{i+1} - x^i \|^2}{\| x^i - x^{i-1} \|} 
 	= \beta \frac{\| x^{i+1} - x^i \|^2}{\| x^i - x^{i-1} \|}
\end{align*}
for all $i\in\{k_0+1,k_0+2,\ldots,k+1\}$,
where we used the constant 
$ \beta := \frac{\delta \gamma_{\min}}{2 (\bar \gamma_{\rho}
+ L_{\rho})} $. Noting that $ a + b \geq 2 \sqrt{ab} $ holds for all real
numbers $ a, b \geq 0 $, we therefore obtain
\begin{equation*}
	\frac{1}{\beta} \Delta_{i,i+1} + \| x^i - x^{i-1} \| 
	\geq 2 \sqrt{\frac{1}{\beta} \Delta_{i,i+1} \| x^i - 
	x^{i-1} \|}
 	\geq 2 \| x^{i+1} - x^i \|	
\end{equation*}
for all $i\in\{k_0+1,k_0+2,\ldots,k+1\}$.
Summation yields
\begin{align*}
	2 \sum_{i=k_0 + 1}^{k+1} \| x^{i+1} - x^i \| 
	& \leq \sum_{i=k_0+1}^{k+1} \| x^i - x^{i-1} \| +
	\frac{1}{\beta} \sum_{i=k_0+1}^{k+1} \Delta_{i,i+1} \\
	& = \sum_{i=k_0+1}^k \| x^{i+1} - x^{i} \| + 
	\| x^{k_0 +1} - x^{k_0} \|   
	+ \frac{1}{\beta} 
	\Delta_{k_0+1, k+2}\\
	& \leq \sum_{i=k_0+1}^{k+1} \| x^{i+1} - x^{i} \| + 
	\| x^{k_0 +1} - x^{k_0} \|   + \frac{1}{\beta} 
	\Delta_{k_0+1, k+2} .
\end{align*}
Subtracting the first summand from the right-hand side, 
exploiting the estimate \eqref{Eq:induktiv}, and using
the nonnegativity as well as monotonicity of the desingularization function 
$ \chi $, we obtain
\[
   \sum_{i=k_0+1}^{k+1} \| x^{i+1} - x^i \| \leq 
   \sqrt{\frac{2 \big( \psi (x^{k_0}) - \psi (x^*) \big)}{\delta
   \gamma_{\min}}} + \frac{1}{\beta} \chi \big( \psi (x^{k_0}) -
   \psi (x^*) \big) .
\]
Adding the term $ \| x^{k_0+1} - x^{k_0} \| $ to both sides 
and using \eqref{Eq:induktiv} once again, we get
\[
   \sum_{i=k_0}^{k+1} \| x^{i+1} - x^i \| \leq 
   \sqrt{\frac{8 \big( \psi (x^{k_0}) - \psi (x^*) \big)}{\delta
   \gamma_{\min}}} + \frac{1}{\beta} \chi \big( \psi (x^{k_0}) -
   \psi (x^*) \big) .
\]
Hence, statement \ref{Item:Ind-2} holds for $ k+1 $ in place of $k$, 
and this completes the induction.

In particular, it follows from \ref{Item:Ind-1} that 
$ x^k \in B_{\alpha} (x^*) $ for all $ k \geq k_0 $. Taking
$ k \to \infty $ in \eqref{Eq:Ind-2} therefore shows that 
$ \{ x^k \} $ is a Cauchy sequence and, thus, convergent.
Since we already know that $ x^* $ is an accumulation point,
it follows that the entire sequence $ \{ x^k \} $ converges
to $ x^* $.
\end{proof}

	Let us note that \cref{Thm:GlobConv} says that, in the
	presence of \cref{Ass:ProxGradMonotone} and the KL property (on the overall domain of $\phi$),
	any sequence $\{x^k\}$ generated by \cref{Alg:MonotoneProxGrad}
	either satisfies $\|x^k\|\to\infty$ or converges to
	a limit point (which is an M-stationary point of \eqref{Eq:P} by \cref{Thm:ConvProxGrad}).
	This alternative behavior, which typically comes along with the KL property,
	see e.g.\ \cite[Theorem~3.2]{AttouchBolteRedontSoubeyran2010},
	has been observed for the first time in
	\cite[Theorem~3.2]{AbsilMahonyAndrews2005} in the context of descent
	methods for analytic functions.

We finally state our rate-of-convergence result for one
particular class of desingularization functions. The 
result holds for a more general class of such functions,
and we comment on this after the proof. To keep the 
notation simple and since this result, having in mind
the previous ones, is more or less a standard observation,
we decided to state this rate-of-convergence result
in the following way.

\begin{theorem}\label{Thm:Rate-of-Conv}
Let \Cref{Ass:ProxGradMonotone} hold, and let $ \{ x^k \} $
be any sequence generated by \Cref{Alg:MonotoneProxGrad}.
Suppose that $ \{ x^k \}_K $ is a subsequence converging to
some limit point $ x^* $, and that $ \psi $ has the KL property
at $ x^* $. Then the entire sequence $ \{ x^k \} $ converges
to $ x^* $, and if the corresponding desingularization function
has the form $ \chi (t) = c t^{1/2} $ for some $ c > 0 $, 
the following statements hold:
\begin{enumerate}
	\item \label{Item:Rate-of-Conv-psi}
	   the sequence $ \{ \psi (x^k ) \} $ converges 
	   $ Q $-linearly to $ \psi (x^*) $,
	\item \label{Item:Rate-of-Conv-x}
	   the sequence $ \{ x^k \} $ converges R-linearly to $ x^* $.
\end{enumerate}
\end{theorem}

\begin{proof}
In view of \Cref{Thm:GlobConv}, we only need to verify the quantitative
statements~\ref{Item:Rate-of-Conv-psi} and~\ref{Item:Rate-of-Conv-x}
of the theorem.

As noted at the beginning of the proof of \cref{Thm:GlobConv}, we may
assume, without loss of generality, that $ \psi (x^k) > \psi (x^*) $
holds for all $ k \in \N $. In view of
\Cref{Lem:psi-conv}, we then have
\begin{equation*}
	x^k \in B_{\alpha} (x^*) \cap \big\{ x \in \dom 
	\phi \mid \psi (x^*) < \psi (x) < \psi (x^*) + \eta \big\}
\end{equation*}
for all $ k \in \N $ sufficiently large, where $ \alpha > 0 $
is the constant from \eqref{Eq:alpha} and $ \eta > 0 $
denotes the constant from the definition of the
desingularization function $\chi$. Since $ \psi $
satisfies the KL property at $ x^* $ with $ \chi (t) = c t^{1/2} $, we have 
\begin{align*}
   1 &\leq \chi ' \big( \psi (x^{k+1}) - \psi (x^*) \big)
   \dist \big( 0, \partial \psi (x^{k+1}) \big)\\
   &=\frac{c}{2} \big( \psi (x^{k+1}) - \psi (x^*) \big)^{-1/2}
   \dist \big( 0, \partial \psi (x^{k+1}) \big)
\end{align*}
for all sufficiently large $ k\in\N $. Taking into account
\Cref{Lem:DistanceSubgrad}, this yields
\[
   1 \leq 
   \frac{c (\bar \gamma_{\rho} + L_{\rho})}{2} 
   \big( \psi (x^{k+1}) - \psi (x^*) \big)^{-1/2}
   \| x^{k+1} - x^k \|
\]
for all $ k \in \N $ sufficiently large, where $ \bar \gamma_{\rho}>0 $
is the constant from \eqref{Eq:gamma-rho} and $L_{\rho}>0$ is
the global Lipschitz constant of $\nabla f$ on $C_\rho$ from \eqref{Eq:C_rho}.
Rearranging this expression gives us 
\begin{equation}\label{Eq:Rate-1}
   \| x^{k+1} - x^k \| \geq \frac{2}{c (\bar \gamma_{\rho} + 
   L_{\rho})} \big( \psi (x^{k+1}) - \psi (x^*) \big)^{1/2} .
\end{equation}
On the other hand, by the acceptance criterion \eqref{Eq:StepCrit}
and $ \gamma_k \geq \gamma_{\min} $, 
we have
\begin{equation}\label{Eq:Rate-2}
	\psi (x^{k+1}) - \psi (x^k) \leq - \delta 
	\frac{\gamma_{\min}}{2} \| x^{k+1} - x^k \|^2 .
\end{equation}
Combining \eqref{Eq:Rate-1} and \eqref{Eq:Rate-2}, we obtain 
\begin{align*}
	\big( \psi (x^{k+1}) - \psi (x^*) \big) - 
	\big( \psi (x^{k}) - \psi (x^*) \big)
	& = \psi (x^{k+1}) - \psi (x^{k}) \\
	& \leq - \delta 
	\frac{\gamma_{\min}}{2} \| x^{k+1} - x^k \|^2 \\
	& \leq -\frac{2 \delta \gamma_{\min}}{c^2 
	(\bar \gamma_{\rho} + L_{\rho})^2}
    \big( \psi (x^{k+1}) - \psi (x^*) \big) \\
    & = - \sigma \big( \psi (x^{k+1}) - \psi (x^*) \big)
\end{align*}
for all $ k\in\N $ sufficiently large, where 
we used the constant 
$ \sigma := \frac{2 \delta \gamma_{\min}}{c^2(\bar \gamma_{\rho} + L_{\rho})^2} $
for brevity. Rearranging these terms,
we find that
\begin{equation}\label{Eq:square-root}
   \psi (x^{k+1}) - \psi (x^*) \leq \frac{1}{1 + \sigma}
   \big( \psi (x^{k}) - \psi (x^*) \big)
\end{equation}
holds for all $ k\in\N$ large enough,
which shows that the sequence $ \{\psi (x^k)\} $ converges 
Q-linearly to $ \psi (x^*) $.

To verify statement \ref{Item:Rate-of-Conv-x}, observe
that the descent test \eqref{Eq:StepCrit} and the monotonicity of
the sequence $ \{ \psi (x^k) \} $ yield
\[
   \frac{\delta \gamma_{\min}}{2} \| x^{k+1} - x^k \|^2 
   \leq \psi (x^k) - \psi (x^{k+1}) \leq \psi (x^k) -
   \psi (x^*) =: \psi_k,
\]
and that the sequence $ \{ \psi_k \} $ is Q-linearly convergent
in view of part \ref{Item:Rate-of-Conv-psi}. Taking this into 
account, it is not difficult to see 
that there exist 
constants $ \omega > 0 $ and $ \mu \in (0,1) $ such that
\[
   \| x^{k+1} - x^k \| \leq \omega \mu^k 
\]
holds for all sufficiently large $ k\in\N $. Hence, for given integers 
$ \ell > k > 0 $ large enough, we therefore obtain
\[
   \| x^{\ell + 1} - x^k \| \leq \sum_{j=k}^{\ell}
   \| x^{j+1} - x^j \| \leq \omega \sum_{j= k}^{\ell} \mu^j 
   \leq \omega \mu^k \sum_{j=0}^{\infty} \mu^j
   = \frac{\omega}{1 - \mu} \mu^k.
\]
Taking the limit $ \ell \to \infty $ yields
\[
   \| x^k - x^* \| \leq \frac{\omega}{1 - \mu} \mu^k
\]
for all large enough $ k\in\N $.  This completes the proof of 
the (local) R-linear convergence of $ \{ x^k \} $ to its
limit $ x^* $.
\end{proof}

We note that similar rate-of-convergence results can be 
obtained for the more general case where the desingularization
function is given by $ \chi (t)  = c t^{\kappa} $ for some
$ \kappa \in (0,1 ] $. The easiest way to see that is to 
modify the previous proof and to apply, for example,
\cite[Lemma~1]{Aragon2018}.

\section{Conclusions}\label{Sec:Final}

In this paper, we have shown that convergence of the whole sequence
generated by proximal gradient methods applied to the composite
optimization problem \eqref{Eq:P} can be achieved whenever
the gradient of the smooth function $f$ is locally Lipschitz
continuous while the objective function $\psi$ possesses the
KL property at all points of its domain. For our analysis, we
neither needed a priori boundedness of iterates and stepsizes
nor any additional convexity assumptions. 
Our findings also gave rise to the statement of associated 
rate-or-convergence results.

In \cite{DeMarchi2022}, the author shows that the global convergence properties
of \cref{Alg:MonotoneProxGrad} from \cref{Thm:ConvProxGrad} remain valid if,
instead of the exploited monotone line search, a nonmonotone scheme is used to
determine the step sizes. In the future, it should be clarified whether the
results of \cref{Thm:GlobConv,Thm:Rate-of-Conv} can be carried over to nonmonotone
proximal gradient methods.

Several generalizations of the proximal gradient method involving, e.g., inertial terms
or Bregman distances, 
see \cite{BauschkeBolteTeboulle2017,BolteSabachTeboulleVaisbourd2018,BotCsetnek2016,BotCsetnekLaszlo2016} 
and the references therein,
have been investigated in the presence of global Lipschitzness of the gradient associated with the smooth
term, as well as the KL property. 
Keeping our findings in mind, it might be promising to check whether our technique of proof
can be applied in these settings to weaken the current Lipschitz assumptions.

\subsection*{Acknowledgements}

The authors wish to thank J\'{e}r\^{o}me Bolte and the reviewers for pointing their attention to several subject-related references. 
 

\end{document}